**УДК 521.98+521.17**

# Эволюционные уравнения для кубических стохастических процессов

Б.Ж. Мамуров

Аннотация

В данной работе по аналогии [1,2] определяется кубический стохастический процесс и изучается эволюция (динамика) некоторой биологической системы $E$ содержащей не менее трех элементов.

Будем считать биологический процесс определенным, если задана тройка $(E,F,\mathrm{M})$, где $F$ есть $\sigma$- алгебра подмножеств $E$, $M$- совокупность вероятностных мер на $(E,F)$ и из знания закона взаимодействия трех элементов $x,y,z$ из E и вероятностную меру $m_{t_0}\in M$ в некоторый момент времени $t_0$ следует знание закона распределения вероятностей $m_{t_1}\in M$ системы $E$ в момент времени $t_1>t_0$. Изучаемый процесс определяется следующим образом: Процесс изменения вероятностных мер во времени назовем кубично стохастически определенным, если при любом выборе элементов $x,y$ и $z$ из E, множества $A\in F$ и моментов времени $t_1$ и $t_2$, причем $t_2\geq t_1+1$, существует определенная вероятность $P(t_1,x,y,z,t_2,A)$ того, что при взаимодействии элементов $x,y$ и $z$ в момент времени $t_1$, в момент времени $t_2$ реализуется один из элементов, принадлежащих множеству $A$ .

При этом, если в момент времени $t_1$ известна вероятностная мера $m_{t_1}\in M$ , то $m_{t_2}(t_2\geq t_1+1)$ определяется по формуле

$$m_{t_2}(A)=\int_E\int_E\int_E P(t_1,x,y,z,t_2,A)m_{t_1}(dx)m_{t_1}(dy)m_{t_1}(dz) \qquad (1)$$

Далее требуется выполнение следующих условий:

(I) $P(t,x,y,z,t+1,A)=P(0,x,y,z,1,A)$ для всех $t\geq 1$;

(II) Значение $P(t_1,x,y,z,t_2,A)$ не меняется при произвольной перестановке аргументов $x,y,z$ для всех $x,y,z\in E$, $A\in F$;

(III) $P(t_1,x,y,z,t_2,A)$- мера на $(E,F)$ для всех фиксированных $x,y,z\in E$.

(IV) $P(t_1,x,y,z,t_2,A)$ как функция трёх переменных $x,y$ и $z$ измерима относительно $(E\times E\times E,\quad F\otimes F\otimes F)$ для всех фиксированных $A\in F$;

(V) Для произвольных $t_1<t_2<t_3$ таких, что $t_2-t_1\geq 1$ и $t_3-t_2\geq 1$ имеет место равенство

$$P(t_1,x,y,z,t_3,A)=\int_E\int_E\int_E P(t_1,x,y,z,t_2)P(t_2,u,\vartheta,q,t_3,A)m_{t_2}(d\vartheta)\cdot m_{t_2}(dq) \qquad (2)$$

Это равенство будем называть фундаментальным уравнением рассматриваемого процесса, оно является аналогом известного уравнения Колмогорова-Чепмена в теории марковских процессов.

Укажем еще одну аналогию с марковскими процессами, мы предложим однородность процесса за нулевой промежуток.

$$P(t,x,y,t,A)=\delta_{x,y}(A)=\begin{cases}1, & если \quad (x,y)\in A\\ 0, & если \quad (x,y)\notin A\end{cases}$$

Условия (1) и требования неравенства $t_2 \ge t_1 + 1$ определения нашего процесса принимаются, исходя из общих положений генетики, которая служит для нас в качестве иллюстрации.

Пусть $E=\{1,2,...,n\}$ и $x^{(0)}=\{x_1^{(0)},x_2^{(0)},\dots,x_n^{(0)}\}$ - начальное распределение на *E*. Введем обозначения:

$P(0,i,j,k,t,l)=p_{ijk,l}$, $\quad p(s,i,j,k,t,l)=p_{ijk,l}^{[s,t]}$

Распределения на момент $t=1$, $x^{(1)}=\left(x_1^{(1)},x_2^{(1)},\dots,x_n^{(1)}\right)$ в силу уравнения (1) определяется следующим образом

$$x_l^{(1)}=\sum_{i,j,k=1}^{n}p_{ijk,l}x_i^{(0)}x_j^{(0)}x_k^{(0)},\quad k=\overline{1,n}.$$

В этом случае условия (I) имеют следующий вид: $p_{ijk,l}^{[t,t+1]}=p_{ijk,l}$.

Предположим что время дискретно. Тогда оператор $p^{[0,2]}$ определяет: $P^{[0,2]}x^{(0)}=Px^{(1)}=x^{(2)}$ и $x_1^{(2)}=\sum_{i,j,k}p_{ijk,l}x_i^{(1)}x_j^{(1)}x_k^{(1)}=\sum_{\alpha,\beta,\theta}p_{\alpha\beta\theta,k}^{[0,2]}x_\alpha^{(0)}x_\beta^{(0)}x_\theta^{(0)}$.

Уравнения (1) примут вид

$$x_l^t=\sum_{i,j,k}p_{ijkl}^{[s,t]}x_i^{(s)}x_j^{(s)}x_k^{(s)} \tag{3}$$

а фундаментальное уравнение (2)

$$p_{ijk,l}^{[s,t]}=\sum_{m,\gamma,\delta}p_{ijk,m}^{[s,\tau]}p_{m\gamma\delta,l}^{[\tau,t]}x_\gamma^{(\tau)}x_\delta^{(\tau)}, \tag{4}$$

где $\tau-s\ge 1$ и $\tau-\varepsilon\ge 1$.

Пример: 1. $E=\{1,2,3\}$ и $x^{(0)}=\left(x_1^{(0)},x_2^{(0)},x_3^{(0)}\right)$-начальное распределение на $E$.

Тогда следующая система вероятностей при $0\le\varepsilon\le\frac{1}{2}$

$$P_{iii,i}^{[s,t]}=\frac{(1-2\varepsilon)^{t-s}}{2^{t-s-1}}\left[(2^{t-s-1}-1)(1-2\varepsilon)^s x_i^{(0)}+1\right],\quad i=1,2,3;$$

$$P_{iii,l}^{[s,t]}=\frac{(2^{t-s-1}-1)}{2^{t-s-1}}x_l^{(0)},\ (i\ne l);$$

а в остальных случаях

$$P_{ijk,l}^{[s,t]}=\frac{1}{2^{t-s-1}}\left[(2^{t-s-1}-1)(1-2\varepsilon)^s x_l^{(0)}+\frac{1}{3}\right]$$

удовлетворяют условиям (I)-(V).

Предположим, что функции $p_{ijk,l}^{[s,t]}$ при $t \geq s+2$ дифференцируемы по $s$ и $t$.

Для $t \geq s+2$ в силу (4) имеем:

$$p_{ijk,l}^{[s,t+\Delta]} - p_{ijk,l}^{[s,t]} = \sum_{m,\gamma,\delta} p_{ijk,m}^{[s,t-1]} p_{m\gamma\delta,l}^{[t-1,t+\Delta]} x_\gamma^{(t-1)} x_\delta^{(t-1)} - \sum_{m,\gamma,\delta} p_{ijk,m}^{[s,t-1]} p_{m\gamma\delta,l}^{[t-1,t]} x_\gamma^{(t-1)} x_\delta^{(t-1)} =$$

$$= \sum_{m\gamma\delta} p_{ijk,m}^{[s,t-1]} \left\{ p_{m\gamma\delta,l}^{[t-1,t+\Delta]} - p_{ijk,m}^{[t-1,t]} \right\} x_\gamma^{(t-1)} x_\delta^{(t-1)} \tag{5}$$

Так как $p_{m\gamma\delta,l}^{[t-1,t]} = p_{m\gamma\delta,l}$ для всех $t$, положим $a_{m\gamma\delta,l}(t) = \lim\limits_{\Delta \to 0} \frac{p_{m\gamma\delta,l}^{[t-1,t+\Delta]} - p_{m\gamma\delta,l}}{\Delta}$ в предположении, что существует правая часть.

Разделив обе стороны равенства (5) на $\Delta$ и переходя к пределу, получим

$$\frac{\partial p_{ijk,l}^{[s,t]}}{\partial t} = \sum_{m,\gamma,\delta} a_{m\gamma\delta,l}(t) x_\gamma^{(t-1)} x_\delta^{(t-1)} p_{ijk,l}^{[s,t-1]} \tag{6}$$

Предположим, что $t > s+2$. В силу (4)

$$p_{ijk,l}^{[s,t]} - p_{ijk,l}^{[s+\Delta,t]} = \sum_{m,\gamma,\delta} \Big[ p_{ijk,m}^{[s,s+1+\Delta]} p_{m\gamma\delta,l}^{[s+1+\Delta,t]} x_\gamma^{(s+1+\Delta)} x_\delta^{(s+1+\Delta)} -$$

$$- p_{ijk,m}^{[s+\Delta,s+\Delta+1]} p_{m\gamma\delta,l}^{[s+\Delta+1,t]} x_\gamma^{(s+\Delta+1)} x_\delta^{(s+\Delta+1)} \Big] = \sum_{m\gamma\delta} \Big[ p_{ijk,m}^{[s,s+\Delta+1]} - p_{ijk,m} \Big] x_\gamma^{(s+\Delta+1)} x_\delta^{(s+\Delta+1)} p_{m\gamma\delta,l}^{[s+\Delta+1,t]}$$

Разделив обе стороны равенства (6) на $\Delta$ и переходя к пределу, получим

$$\frac{\partial p_{ijk,l}^{[s,t]}}{\partial s} = -\sum_{m\gamma\delta} a_{m\gamma\delta,l}(s+1) x_\gamma^{(s+1)} x_\delta^{(s+1)} P_{m\gamma\delta,l}^{[s+1,t]} \tag{7}$$

Теперь рассмотрим дифференциальные уравнения относительно распределений $x_k^{(t)}$. Так как в силу (3) $x_k^{(t+\Delta)} = \sum_{i,j,l} p_{ijl,k}^{[t-1,t+\Delta]} x_i^{(t-1)} x_j^{(t-1)} x_l^{(t-1)}$ и $x_k^{(t)} = \sum_{ij,l} p_{ijl,k}^{[t-1,t]} x_i^{(t-1)} x_j^{(t-1)} x_l^{(t-1)}$.

Для $t \geq 1$ имеем

$$x_k^{(t+1)} - x_k^{(t)} = \sum_{ij,l} \Big[ p_{ijl,k}^{[t-1,t+\Delta]} - p_{ijl,k}^{[t-1,t]} \Big] x_i^{(t-1)} x_j^{(t-1)} x_l^{(t-1)} \tag{8}$$

Разделив обе стороны равенства (8) на $\Delta$ и переходя к пределу, получим

$$\frac{dx_k^{(t)}}{dt} = \sum_{ij,l} a_{ijl,k}(t) x_i^{(t-1)} x_j^{(t-1)} x_l^{(t-1)}. \tag{9}$$

Система дифференциальных уравнений (6) и (9) являются уравнениями с запаздывающим аргументом, а система дифференциальных уравнений (7)- уравнением с опережающим аргументом.

Пусть $E$-континуальна. Рассмотрим кубично стохастически определенный процесс, заданный на измеримом пространстве $(E, F)$ с начальной мерой $m_0$.

При $t > s+2$ в силу (V) имеет

$$P(s,x,y,z,t+\Delta,A)-P(s,x,y,z,t,A)=$$
$$=\iiint_{E\,E\,E} P(s,x,y,z,t-1,du)\{P(t-1,u,\vartheta,q,t+\Delta,A)-P(t-1,u,\vartheta,q,t,A)\}m_{t-1}(d\vartheta)m_{t-1}(dq)$$

Предположим, что $C(t,u,\vartheta,q,A)=\lim\limits_{\Delta\to 0}\dfrac{P(t-1,u,\vartheta,q,t+\Delta,A)-P(t-1,u,\vartheta,q,t,A)}{\Delta}$ существует.

Тогда переходя к пределу при $\Delta\to 0$

$$\frac{\partial P(s,x,y,z,t,A)}{\partial t}=\iiint_{E\,E\,E} P(s,x,y,z,t-1)du\ C(t,u,\vartheta,q,A)m_{t-1}(d\vartheta)m_{t-1}(dq) \qquad (10)$$

получим первое интегро-дифференциальное уравнение.

Аналогично, из $P(s,x,y,z,t,A)-P(s+\Delta,x,y,z,t,A)=\iiint_{E\,E\,E}\{P(s,x,y,z,s+1+\Delta,du)-$ $-P(s,x,y,z,s+1,du)\}P(s+1+\Delta,u,\vartheta,q,t,A)m_{s+1+\Delta}(d\vartheta)m_{s+1+\Delta}(dq)$ получим второй интегро-дифференциальное уравнение

$$\frac{\partial P(s,x,y,z,t,A)}{\partial s}=-\iiint_{E\,E\,E} C(s+1,x,y,z,du)P(s+1,u,\vartheta,q,t,A)m_{s+1}(d\vartheta)m_{s+1}(dq). \qquad (11)$$

Пусть $E=R$ и $A_w=(-\infty;w]$ для любого $w\in R$. Положим $F(s,x,y,z,t,w)=P(s,x,y,z,t,A_w)$.

Видно что $F(s,x,y,z,t,w)$ как функция от $w$ монотонна, непрерывна слева и удовлетворяет граничным условиям $F(s,x,y,z,t,-\infty)=0$, $F(s,x,y,z,t,+\infty)=1$. Для функции $F(s,x,y,z,t,w)$ условие (V) будет имеет вид:

$$F(s,x,y,z,t,w)=\iiint_{E\,E\,E} dF(s,x,y,z,\tau,u)F(\tau,u,\vartheta,q,t,w)m_\tau(d\vartheta)m_\tau(dq).$$

Если $F(s,x,y,z,t,w)$ как функция $w$ абсолютно непрерывна, то $F(s,x,y,z,t,w)=\int_{-\infty}^{w} f(s,x,y,z,t,u)du$. Где $f(s,x,y,z,t,w)$-неотрицательная функция измерима относительно аргументов $x,y,z,w$ и удовлетворяет условиям

$$\int_{-\infty}^{\infty} f(s,x,y,z,t,w)dw=1 \quad \text{и}$$

$$f(s,x,y,z,t,w)=\iiint_{E\,E\,E} f(s,x,y,z,\tau,u)f(\tau,u,\vartheta,q,t,w)m_\tau(d\vartheta)m_\tau(dq)du \qquad (12)$$

(10) и (11) для $F(s,x,y,z,t,w)$ имеет следующий вид:

$$\frac{\partial F(s,x,y,z,t,w)}{\partial t}=\iiint_{E\,E\,E}\frac{\partial F(s,x,y,z,t-1,u)}{\partial u}C(t,u,\vartheta,q,w)du\,m_{t-1}(d\vartheta)m_{t-1}(dq)$$

$$\frac{\partial F(s,x,y,z,t,w)}{\partial s}=\iiint_{E\,E\,E}\frac{\partial C(s+1,x,y,z,u)}{\partial u}F(s+1,u,\vartheta,q,t,w)du\,m_{s+1}(d\vartheta)m_{s+1}(dq).$$

Если $\overline{C}(s+1,x,y,z,u)=\dfrac{\partial C(s+1,x,y,z,u)}{\partial u}$ тогда,

$$\frac{\partial F(s,x,y,z,t,w)}{\partial s}=\int\limits_E\int\limits_E\int\limits_E \overline{C}(s+1,x,y,z,u)F(s+1,u,\vartheta,q,t,w)dum_{s+1}(d\vartheta)m_{s+1}(dq).$$

Если $a(t,u,\vartheta,q,w)=\lim\limits_{\Delta\to 0}\dfrac{f(t-1,u,\vartheta,q,t+\Delta,w)-f(0,u,\vartheta,q,1,w)}{\Delta}$ предел существует, тогда получим

$$\frac{\partial f(s,x,y,z,t,w)}{\partial t}=\int\limits_E\int\limits_E\int\limits_E a(t,u,\vartheta,q,w)f(s,x,y,z,t-1,u)dum_{t-1}(d\vartheta)m_{t-1}(dq) \quad (13)$$

$$\frac{\partial f(s,x,y,z,t,w)}{\partial t}=-\int\limits_E\int\limits_E\int\limits_E a(s+1,x,y,z,u)f(s+1,u,\vartheta,q,t,z)dum_{s+1}(d\vartheta)m_{s+1}(dq)$$

(14).

При некоторых условиях интегро-дифференциальные уравнения (13) и (14) сводятся к дифференциальным. Пусть $t>s+2$ , тогда из (12) имеем

$$f(s,x,y,z,t,w)-f(s+\Delta,x,y,z,t,w)=\int\limits_E\int\limits_E\int\limits_E\{f(s,x,y,z,s+1+\Delta,u)-$$
$$-f(s+\Delta,x,y,z,s+1+\Delta,u)\}f(s+1+\Delta,u,\vartheta,q,t,w)m_{s+1+\Delta}(d\vartheta)m_{s+1+\Delta}(dq)du$$

(15)

Разложив функцию $f(s+1+\Delta,u,\vartheta,q,t,w)$ в ряд Тейлора в окрестности точки $(x,y,z)$, имеем:

$$f(s+1+\Delta,u,\vartheta,q,t,w)=f(s+1+\Delta,x,y,z,t,w)+\frac{\partial f(s+1+\Delta,x,y,z,t,w)}{\partial u}(u-x)+$$
$$+\frac{\partial f(s+1+\Delta,x,y,z,t,w)}{\partial\vartheta}(\vartheta-y)+\frac{\partial f(s+1+\Delta,x,y,z,t,w)}{\partial q}(q-z)+$$
$$+\frac{1}{2}\frac{\partial^2 f(s+1+\Delta,x,y,z,t,w)}{\partial u^2}(u-x)^2+\frac{1}{2}\frac{\partial^2 f(s+1+\Delta,x,y,z,t,w)}{\partial\vartheta^2}(\vartheta-y)^2+$$
$$+\frac{1}{2}\frac{\partial^2 f(s+1+\Delta,x,y,z,t,w)}{\partial q^2}(q-z)^2+\frac{\partial^2 f(s+1+\Delta,x,y,z,t,w)}{\partial u\partial\vartheta}(u-x)(\vartheta-y)+$$
$$+\frac{\partial^2 f(s+1+\Delta,x,y,z,t,w)}{\partial u\partial q}(u-x)(q-z)+\frac{\partial^2 f(s+1+\Delta,x,y,z,t,w)}{\partial\vartheta\partial q}(\vartheta-y)(q-z)+$$
$$+\frac{1}{2}\frac{\partial^3 f(s+1+\Delta,x,y,z,t,w)}{\partial u^2\partial\vartheta}(u-x)^2(\vartheta-y)+\frac{1}{2}\frac{\partial^3 f(s+1+\Delta,x,y,z,t,w)}{\partial u^2\partial q}(u-x)^2(q-z)+$$

$$+\frac{1}{2}\frac{\partial^3 f(s+1+\Delta,x,y,z,t,w)}{\partial\vartheta^2\partial u}(\vartheta-y)^2(u-x)+\frac{1}{2}\frac{\partial^3 f(s+1+\Delta,x,y,z,t,w)}{\partial\vartheta^2\partial q}(\vartheta-y)^2(q-z)+$$

$$+\frac{1}{2}\frac{\partial^3 f(s+1+\Delta,x,y,z,t,w)}{\partial q^2\partial u}(q-z)^2(u-x)+\frac{1}{2}\frac{\partial^3 f(s+1+\Delta,x,y,z,t,w)}{\partial q^2\partial\vartheta}(q-z)^2(\vartheta-y)+$$

$$+\frac{1}{6}\frac{\partial^3 f(s+1+\Delta,x,y,z,t,w)}{\partial q^3}(u-x)^3+\frac{1}{6}\frac{\partial^3 f(s+1+\Delta,x,y,z,t,w)}{\partial\vartheta^2}(\vartheta-y)^3+$$

$$+\frac{1}{6}\frac{\partial^3 f(s+1+\Delta,x,y,z,t,w)}{\partial q^2}(q-z)^3+\frac{\partial^3 f(s+1+\Delta,x,y,z,t,w)}{\partial u\partial\vartheta\partial q}(u-x)(\vartheta-y)(q-z)$$

16)

Подставим (16) в (15) и рассмотрим те слагаемых которые отличны от нуля:

$$\int_E\int_E\int_E\{f(s,x,y,z,s+1+\Delta,u)-f(s+\Delta,x,y,z,s+1+\Delta,u)\}\frac{\partial f(s+1+\Delta,x,y,z,t,w)}{\partial x}(u-x)\times$$
$$\times m_{s+1+\Delta}(d\vartheta)m_{s+1+\Delta}(dq)du=\frac{\partial f(s+1+\Delta,x,y,z,t,w)}{\partial x}\int_E\{f(s,x,y,z,s+1+\Delta,u)-$$
$$-f(s+\Delta,x,y,z,s+1+\Delta,u)\}(u-x)du$$

Положим

$$a(s,x,y,z,\Delta)=\int_E\{f(s,x,y,z,s+1+\Delta,u)-f(s+\Delta,x,y,z,s+1+\Delta,u)\}(u-x)du\,.$$

Теперь перейдем к слагаемым со вторыми производными

$$\int_E\int_E\int_E\{f(s,x,y,z,s+1+\Delta,u)-f(s+\Delta,x,y,z,s+1+\Delta,u)\}\frac{1}{2}\frac{\partial^2 f(s+1+\Delta,x,y,z,w)}{\partial x^2}\times$$
$$(u-x)^2 m_{s+1+\Delta}(d\vartheta)m_{s+1+\Delta}(dq)du=\frac{1}{2}\frac{\partial^2 f(s+1+\Delta,x,y,z,w)}{\partial x^2}\int_E\{f(s,x,y,z,s+1+\Delta,u)-$$
$$-f(s+\Delta,x,y,z,s+1+\Delta,u)\}(u-x)^2\,du\,.$$

Положим

$$b^2(s,x,y,z,\Delta)=\int_E\{f(s,x,y,z,s+1+\Delta,u)-f(s+\Delta,x,y,z,s+1+\Delta,u)\}(u-x)^2\,du$$

$$\int_E\int_E\int_E\{f(s,x,y,z,s+1+\Delta,u)-f(s+\Delta,x,y,z,s+1+\Delta,u)\}\frac{\partial^2 f(s+1+\Delta,x,y,z,w)}{\partial x\partial y}(u-x)(\vartheta-y)\times$$
$$m_{s+1+\Delta}(d\vartheta)m_{s+1+\Delta}(dq)du=\frac{\partial^2 f(s+1+\Delta,x,y,z,t,w)}{\partial x\partial y}\int\{f(s,x,y,z,s+1+\Delta,u)-$$
$$-f(s+\Delta,x,y,z,s+1+\Delta,u)\}(u-x)du\int_E(\vartheta-y)m_{s+1+\Delta}(d\vartheta)\,.$$

Положим $\int_E(\vartheta-y)m_{s+1+\Delta}(d\vartheta)=\alpha(s+1,y,\Delta)$. Далее

$$\int_E\int_E\int_E\{f(s,x,y,z,s+1+\Delta,u)-f(s+\Delta,x,y,z,s+1+\Delta,u)\}\frac{\partial^2 f(s+1+\Delta,x,y,z,t,w)}{\partial x\partial z}(u-x)\times$$

$$(q-z)m_{s+1+\Delta}(d\vartheta)m_{s+1+\Delta}(dq)du=\frac{\partial^2 f(s+1+\Delta,x,y,z,t,w)}{\partial x\partial z}\int_E\{f(s,x,y,z,s+1+\Delta,u)-$$
$$-f(s+\Delta,x,y,z,s+1+\Delta,u)\}\int_E(u-x)du\int_E(q-z)m_{s+1+\Delta}(dq)\,.$$

$$\int_E\int_E\int_E \{f(s,x,y,z,s+1+\Delta,u)-f(s+\Delta,x,y,z,s+1+\Delta,u)\}\frac{1}{2}\frac{\partial^3 f(s+1+\Delta,x,y,z,t,w)}{\partial x^2\partial y}(u-x)^2\times$$

$$(\vartheta-y)m_{s+1+\Delta}(d\vartheta)m_{s+1+\Delta}(dq)du=\frac{1}{2}\frac{\partial^3 f(s+1+\Delta,x,y,z,t,w)}{\partial x^2\partial y}\int_E\{f(s,x,y,z,s+1+\Delta,u)-$$

$$-f(s+\Delta,x,y,z,s+1+\Delta,u)\}(u-x)^2du\int(\vartheta-y)m_{s+1+\Delta}(d\vartheta)=\frac{1}{2}\frac{\partial^3 f(s+1+\Delta,x,y,z,t,w)}{\partial x^2\partial y}\times$$

$b^2(s,x,y,z,\Delta)\alpha(s+1,y,\Delta)$.

$$\int_E\int_E\int_E \{f(s,x,y,z,s+1+\Delta,u)-f(s+\Delta,x,y,z,s+1+\Delta,u)\}\frac{1}{2}\frac{\partial^3 f(s+1+\Delta,x,y,z,t,w)}{\partial x^2\partial z}(u-x)^2\times$$

$$(q-z)m_{s+1+\Delta}(d\vartheta)m_{s+1+\Delta}(dq)du=\frac{1}{2}\frac{\partial^3 f(s+1+\Delta,x,y,z,t,w)}{\partial x^2\partial z}\int_E\{f(s,x,y,z,s+1+\Delta)-$$

$$-f(s+\Delta,x,y,z,s+1+\Delta,u)\}(u-x)^2du\int(q-z)m_{s+1+\Delta}(dq)=\frac{1}{2}\frac{\partial^3 f(s+1+\Delta,x,y,z,t,w)}{\partial x^2\partial z}\times$$

$b^2(s,x,y,z,\Delta)\alpha(s+1,z,\Delta)$.

$$\int_E\int_E\int_E \{f(s,x,y,z,s+1+\Delta,u)-f(s+\Delta,x,y,z,s+1+\Delta,u)\}\frac{1}{2}\frac{\partial^3 f(s+1+\Delta,x,y,z,t,w)}{\partial^2 y\partial x}(\vartheta-y)^2\times$$

$$(u-x)m_{s+1+\Delta}(d\vartheta)m_{s+1+\Delta}(dq)du=\frac{1}{2}\frac{\partial^3 f(s+1+\Delta,x,y,z,t,w)}{\partial^2 y\partial x}\int\{f(s,x,y,z,s+1+\Delta,u)-$$

$$-f(s+\Delta,x,y,z,s+1+\Delta,u)\}(u-x)du\int(\vartheta-y)^2m_{s+1+\Delta}(d\vartheta)=\frac{1}{2}a(s,x,y,z,\Delta)\times$$

$\dfrac{\partial^3 f(s+1+\Delta,x,y,z,t,w)}{\partial^2 y\partial x}\cdot\alpha_2(s+1,y,\Delta)$ где $\alpha_2(s+1,y,\Delta)=\int(\vartheta-y)^2m_{s+1+\Delta}(d\vartheta)$.

$$\int_E\int_E\int_E \{f(s,x,y,z,s+1+\Delta,u)-f(s+\Delta,x,y,z,s+1+\Delta,u)\}\frac{1}{2}\frac{\partial^3 f(s+1+\Delta,x,y,z,t,w)}{\partial z^2\partial x}(q-z)^2\times$$

$$(u-x)m_{s+1+\Delta}(d\vartheta)m_{s+1+\Delta}(dq)du=\frac{1}{2}\frac{\partial^3 f(s+1+\Delta,x,y,z,t,w)}{\partial z^2\partial x}\int\{f(s,x,y,z,s+1+\Delta,u)-$$

$$-f(s+\Delta,x,y,z,s+1+\Delta,u)\}(u-x)du\int_E(q-z)^2m_{s+1+\Delta}(d\vartheta)=\frac{1}{2}\frac{\partial^3 f(s+1+\Delta,x,y,z,t,w)}{\partial z^2\partial x}\times$$

$a(s,x,y,z,\Delta)\alpha_2(s+1,z,\Delta)$.

$$\int_E\int_E\int_E \{f(s,x,y,z,s+1+\Delta,u)-f(s+\Delta,x,y,z,s+1+\Delta,u)\}\frac{\partial^3 f(s+1+\Delta,x,y,z,t,w)}{\partial z\partial y\partial z}(u-x)\times$$

$$(\vartheta - y)(q-z)m_{s+1+\Delta}(d\vartheta)m_{s+1+\Delta}(dq)du = \frac{\partial^3 f(s+1+\Delta,x,y,z,t,w)}{\partial x\partial y\partial z}\int\{f(s,x,y,z,s+1+\Delta,u) -$$

$$-f(s+\Delta,x,y,z,s+1+\Delta,u)\}(u-x)du\int_E(\vartheta - y)m_{s+1+\Delta}(d\vartheta)\cdot\int_E(q-z)m_{s+1+\Delta}(dq) =$$

$$= \frac{\partial^3 f(s+1+\Delta,x,y,z,t,w)}{\partial x\partial y\partial z}a(s,x,y,z,\Delta)\alpha(s+1,y,\Delta)\alpha(s+1,z,\Delta).$$

Так как остальные слагаемые равны нулю:

$$f(s,x,y,z,t,w) - f(s+\Delta,x,y,z,t,w) = a(s,x,y,z,\Delta)\frac{\partial f(s+1+\Delta,x,y,z,t,w)}{\partial x} +$$

$$+b^2(s,x,y,z,\Delta)\cdot\frac{1}{2}\frac{\partial^2 f(s+1+\Delta,x,y,z,t,w)}{\partial x^2} + \frac{1}{2}a(s,x,y,z,\Delta)\alpha(s+1,y,\Delta)\times$$

$$\frac{\partial^2 f(s+1+\Delta,x,y,z,t,w)}{\partial x\partial y} + a(s,x,y,z,\Delta)\frac{\partial^2 f(s+1+\Delta,x,y,z,t,w)}{\partial x\partial z}\alpha(s+1,z,\Delta) +$$

$$+\frac{1}{2}\frac{\partial^3 f(s+1+\Delta,x,y,z,t,w)}{\partial x^2\partial y}b^2(s,x,y,z,\Delta)\alpha(s+1,y,\Delta) + \frac{1}{2}\frac{\partial^3 f(s+1+\Delta,x,y,z,t,w)}{\partial x^2\partial z}\times$$

$$b^2(s,x,y,\Delta)\alpha(s+1,z,\Delta) + \frac{1}{2}a(s,x,y,z,\Delta)\frac{\partial^3 f(s+1+\Delta,x,y,z,t,w)}{\partial y^2\partial x}\times$$

$$\times\alpha_2(s+1,y,\Delta) + +\frac{1}{2}\frac{\partial^3 f(s+1+\Delta,x,y,z,t,w)}{\partial z^2\partial x}a(s,x,y,z,\Delta)\alpha_2(s+1,z,\Delta) + a(s,x,y,z,\Delta)\times$$

$$\times\alpha(s+1,y,\Delta)\alpha(s+1,z,\Delta)\frac{\partial^3 f(s+1+\Delta,x,y,z,t,w)}{\partial x\partial y\partial z}. \qquad (17)$$

Обе стороны последнего равенства будем разделим на $\Delta$ и переходя к пределу в (17) при $\Delta \to 0$, положим

$$A(s,x,y,z) = \lim_{\Delta\to 0}\frac{a(s,x,y,\Delta)}{\Delta}, \qquad B^2(s,x,y,z) = \lim_{\Delta\to 0}\frac{b^2(s,x,y,\Delta)}{2\Delta},$$

$$D(s+1,y) = \lim_{\Delta\to 0}\alpha(s+1,y,\Delta), \qquad D_2(s+1,y) = \lim_{\Delta\to 0}\frac{\alpha_2(s+1,y,\Delta)}{2\Delta}$$

предполагаем, что эти пределы существуют.

Тогда из (17) имеем

$$\frac{\partial f(s,x,y,z,t,w)}{\partial s} = -A(s,x,y,z)\frac{\partial f(s+1,x,y,z,t,w)}{\partial x} - B^2(s,x,y,z)\frac{\partial^2 f(s+1,x,y,z,t,w)}{\partial x^2} -$$

$$-\frac{1}{2}A(s,x,y,z)D(s+1,y)\frac{\partial^2 f(s+1,x,y,z,t,w)}{\partial x\partial y} - A(s,x,y,z)D(s+1,z)\frac{\partial^2 f(s+1,x,y,z,t,w)}{\partial x\partial y} -$$

$$-A(s,x,y,z)D(s+1,z)\frac{\partial^2 f(s+1,x,y,z,t,w)}{\partial x\partial z} - B^2(s,x,y,z)D(s+1,y)\frac{\partial^3 f(s+1,x,y,z,t,w)}{\partial x^2\partial z} -$$

$$-B^2(s,x,y,z)D(s+1,z)\frac{\partial^3 f(s+1,x,y,z,t,w)}{\partial x^2\partial z}-\frac{1}{2}A(s,x,y,z)D_2(s+1,y)\frac{\partial^3 f(s+1,x,y,z,t,w)}{\partial y^2\partial x}-$$

$$-\frac{1}{2}A(s,x,y,z)D_2(s+1,z)\frac{\partial^3 f(s+1,x,y,z,t,w)}{\partial z^2\partial x}-A(s,x,y,z)D(s+1,y)D(s+1,z)\times$$

$$\times\frac{\partial^3 f(s+1,x,y,z,t,w)}{\partial x\partial y\partial z}. \qquad (18)$$

Из уравнения (12), имеем

$$f(s,x,y,z,t+\Delta,w)-f(s,x,y,z,t,w)=\int\limits_E\int\limits_E\int\limits_E f(s,x,y,z,t-1,u)(f(t-1,u,\vartheta,q,t+\Delta,w)-$$

$$-f(t-1,u,\vartheta,q,t,w))m_{t-1}(d\vartheta)m_{t-1}(dq)du \qquad (19)$$

Пусть для функции $f(s,x,y,z,t-1,w)$ существуют частные производные до 3-го порядка. Разложим ее в ряд Тейлора в окрестности в точки $w$:

$$f(s,x,y,z,t-1,w)=f(s,x,y,z,t-1,w)+\frac{\partial f(s,x,y,z,t-1,w)}{\partial w}(u-w)+\frac{\partial^2 f(s,x,y,z,t-1,w)}{\partial w^2}\times$$

$$\times\frac{(u-w)^2}{2}+\theta\frac{(u-w)^3}{6}.$$

Тогда из (19) получим

$$\frac{f(s,x,y,z,t+\Delta,w)-f(s,x,y,z,t,w)}{\Delta}=f(s,x,y,z,t-1,w)\times$$

$$\times\int\limits_E\int\limits_E\int\limits_E\frac{f(t-1,u,\vartheta,q,t+\Delta,w)-f(t-1,u,\vartheta,q,t,w)}{\Delta}m_{t-1}(d\vartheta)m_{t-1}(dq)du+$$

$$+\frac{\partial f(s,x,y,z,t-1,w)}{\partial w}\int\limits_E\int\limits_E\int\limits_E\frac{f(t-1,u,\vartheta,q,t+\Delta,w)-f(t-1,u,\vartheta,q,t,w)}{\Delta}(u-w)m_{t-1}(d\vartheta)m_{t-1}(dq)du+$$

$$+\frac{\partial^2 f(s,x,y,z,t-1,w)}{\partial w^2}\int\limits_E\int\limits_E\int\limits_E\frac{f(t-1,u,\vartheta,q,t+\Delta,w)-f(t-1,u,\vartheta,q,t,w)}{\Delta}(u-w)^2m_{t-1}(d\vartheta)m_{t-1}(dq)du+$$

$$+\theta\frac{\int\limits_E\int\limits_E\int\limits_E f(t-1,u,v,q,t+\Delta,w)-f(t-1,u,v,q,t,w)}{6\Delta}(u-w)^3m_{t-1}(d\vartheta)m_{t-1}(dq)du$$

$$(20)$$

Предположим, что существуют следующие пределы.

$$\lim_{\Delta\to 0}\frac{\int\limits_E\int\limits_E\int\limits_E f(t-1,u,\vartheta,q,t+\Delta,w)-f(t-1,u,\vartheta,q,t,w)}{\Delta}m_{t-1}(d\vartheta)m_{t-1}(dq)du=\widetilde{N}(t,w)$$

$$\lim_{\Delta\to 0}\frac{\int\limits_E\int\limits_E\int\limits_E\{f(t-1,u,\vartheta,q,t+\Delta,w)-f(t-1,u,\vartheta,q,t,w)\}m_{t-1}(d\vartheta)m_{t-1}(dq)du}{\Delta}=\widetilde{A}(t,w)$$

$$\lim_{\Delta\to 0}\frac{\int\limits_E\int\limits_E\int\limits_E\{f(t-1,u,\vartheta,q,t+\Delta,w)-f(t-1,u,\vartheta,q,t,w)\}(u-w)^2m_{t-1}(d\vartheta)m_{t-1}(dq)du}{2\Delta}=\widetilde{B}^2(t,w)$$

$$\lim_{\Delta\to 0}\frac{\int\limits_E\int\limits_E\int\limits_E\{f(t-1,u,\vartheta,q,t+\Delta,w)-f(t-1,u,\vartheta,q,t,w)\}|u-w|^3 m_{t-1}(d\vartheta)m_{t-1}(dq)du}{\Delta}=0$$

Тогда переходя к пределу в (20) при $\Delta\to 0$ получим дифференциальное уравнение с отклоняющимся аргументом

$$\frac{\partial f(s,x,y,z,t,w)}{\partial w}=\widetilde{N}(t,w)f(s,x,y,z,t-1,w)+\widetilde{A}(t,w)\frac{\partial f(s,x,y,z,t-1,w)}{\partial z}+$$

$$+\widetilde{B}^2(t,w)\frac{\partial f(s,x,y,z,t-1,w)}{\partial w^2}.$$

Пример 2. Семейство функций

$$f(s,x,y,z,t,w)=\frac{\exp\left(-\frac{w-x-y-z}{2^{t+1}-2^{s+2}+1}\right)}{\sqrt{(2^{t+1}-2^{s+2}+1)\pi}}$$

определяет кубический стохастический процесс.

**Литература**.